\documentclass[smallextended,envcountsect]{svjour3}
\usepackage{graphicx,color,cite}
\usepackage{mathrsfs}
\usepackage{amsmath, amscd, amsfonts, amssymb, graphicx, color}
\usepackage[bookmarksnumbered, colorlinks, plainpages]{hyperref}
\usepackage{mathptmx}   
\usepackage[colorinlistoftodos]{todonotes}
\usepackage[left=4cm, right=4cm, top=3.5cm, bottom=3.5cm]{geometry}

\usepackage{latexsym}

\if{
\newtheorem{theorem}{Theorem}[section]

\newtheorem{remark}{Remark}[section]
\newtheorem{example}{Example}[section]

\newtheorem{corollary}{Corollary}[section]
\newtheorem{proposition}{Proposition}[section]

}\fi

\renewcommand\theenumi{\roman{enumi}}

\newcounter{mycount}

\smartqed

\usepackage{etex}

\usepackage{mathtools}

\usepackage{enumitem}
\renewcommand\theenumi{(\roman{enumi})}
\renewcommand\theenumii{(\alph{enumii})}

\usepackage{csquotes}

\usepackage{todonotes}

\makeatletter
\let\orgdescriptionlabel\descriptionlabel
\renewcommand*{\descriptionlabel}[1]{
 \let\orglabel\label
 \let\label\@gobble
 \phantomsection
 \edef\@currentlabel{#1}
 \let\label\orglabel
 \orgdescriptionlabel{#1}
}
\def\th@plain{
 \thm@notefont{}
 \itshape
}
\def\th@definition{
 \thm@notefont{}
 \normalfont
}

\g@addto@macro\th@definition{\thm@headpunct{}}
\g@addto@macro\th@plain{\thm@headpunct{}}
\makeatother

\usepackage{subfig}

\usepackage[final]{showkeys}

\usepackage{etoolbox}

\usepackage{mathrsfs}

\usepackage[titletoc]{appendix}
\usepackage[doc]{optional}

\usepackage{soul}

\usepackage{colortbl,booktabs,
multirow}
\usepackage{xcolor}

\usepackage{cancel}

\usepackage{empheq}

\definecolor{myblue}{rgb}{.8, .8, 1}

\usepackage{hyperref}
\hypersetup{
colorlinks=true,
linkcolor=blue,
citecolor=blue,
filecolor=magenta,
urlcolor=cyan
}

\usepackage[
open,
openlevel=2,
atend,
numbered
]{bookmark}

\usepackage[capitalize, nameinlink, noabbrev]{cleveref}
\crefname{equation}{}{}
\crefname{chapter}{Chapter}{Chapters}
\crefname{item}{item}{items}
\crefname{figure}{Figure}{Figures}
\crefname{theorem}{Theorem}{Theorems}
\crefname{lemma}{Lemma}{Lemmas}
\crefname{proposition}{Proposition}{Propositions}
\crefname{corollary}{Corollary}{Corollarys}
\crefname{definition}{Definition}{Definitions}
\crefname{fact}{Fact}{Facts}
\crefname{example}{Example}{Examples}
\crefname{algorithm}{Algorithm}{Algorithms}
\crefname{remark}{Remark}{Remarks}
\crefname{note}{Note}{Notes}
\crefname{notation}{Notation}{Notations}
\crefname{case}{Case}{Cases}
\crefname{exercise}{Exercise}{Exercises}
\crefname{question}{Question}{Questions}
\crefname{claim}{Claim}{Claims}
\crefname{enumi}{}{}

\usepackage{pgf}

\usepackage{array}
\usepackage{tabu}

\allowdisplaybreaks

\numberwithin{equation}{section}

\renewcommand\theenumi{(\roman{enumi})}
\renewcommand\theenumii{(\alph{enumii})}

\spnewtheorem*{Proof}{Proof.}{\bf}{\rm}
\begin{document}

\title{Primal Characterizations of Stability of Error Bounds for Semi-infinite Convex Constraint Systems in Banach Spaces \thanks{Research  of the first author was supported by  the National Natural Science Foundations of  China (Grant Nos. 11971422 and 12171419), and funded by Science and Technology Project of Hebei Education Department (No. ZD2022037) and the Natural Science Foundation of Hebei Province (A2022201002). Research of the second author benefited from the support of the FMJH Program PGMO and from the support of EDF.}}

\titlerunning{Primal Characterizations of Stability of Error Bounds of Semi-infinite Convex Constraint Systems}

\author{Zhou Wei  \and Michel Th\'era \and Jen-Chih Yao}

\institute{Zhou Wei\at Hebei Key Laboratory of Machine Learning and Computational Intelligence \& College of Mathematics and Information Science, Hebei University, Baoding, 071002, China\\\email{weizhou@hbu.edu.cn}\\
Michel Th\'era \at ORCID 0000-0001-9022-6406 \at XLIM UMR-CNRS 7252, Universit\'e de Limoges, Limoges, France  \\ and \at Federation University Australia, Ballarat\\\email{michel.thera@unilim.fr}\\
Jen-Chih Yao \at Research Center for Interneural Computing, China Medical University Hospital,
China Medical University, Taichung, Taiwan
\\ \email{yaojc@mail.cmu.edu.tw}
}

\date{Received: date / Accepted: date}

\maketitle

\begin{abstract}

This article is devoted to the stability of error bounds (local and global) for semi-infinite convex constraint systems in   Banach spaces. We provide primal characterizations of the stability of local and global error bounds  when systems are subject to  small perturbations.  These characterizations are given in terms of  the directional derivatives of the functions that enter into the definition of these systems. It is proved that the stability of   error bounds is essentially equivalent to verifying that the optimal values of several minimax problems, defined in terms of the directional derivatives of the functions defining  these systems, are outside of some neighborhood of zero. Moreover, such stability only requires that all component functions entering the system have the same linear perturbation. When these stability results are applied to the sensitivity analysis of Hoffman's constants for semi-infinite linear systems, primal criteria for Hoffman's constants to be uniformly bounded under perturbations of the problem data are obtained.


\end{abstract}

\section{Introduction}

This paper is devoted to the study of  the  stability of error bounds for semi-infinite convex systems under data perturbations in a Banach space $\mathbb{X}$. By a semi-infinite convex system, we mean   the problem of finding $x\in\mathbb{X}$ satisfying
\begin{equation}\label{1}
  f_t(x)\leq 0 \ \ {\rm for\ all} \ \ t\in T,
\end{equation}
where $T$ is a compact, possibly infinite Hausdorff space, $f_t:\mathbb{X}\rightarrow\mathbb{R}, t\in T$ are convex functions such that $(t,x)\mapsto f_t(x)$ is continuous on $T\times\mathbb{X}$. The aim   of this work    is to investigate the behavior of error bounds  of system \eqref{1} under  perturbations of its  data.

Research on error bounds can be traced back to Hoffman's pioneering work in \cite{28}. More precisely,   given a $m\times n$ matrix $A$ and a  $m$-vector $b$,  in order  to estimate the Euclidean distance from a
point $x$ to its nearest point in the polyhedral set $\{u:Au\leq b\}$,  Hoffman \cite{28} showed that this distance is bounded from above by some constant (only dependent on $A$) times the Euclidean norm of the residual error $\|(Ax-b)_+\|$, where $(Ax-b)_+$ denotes the positive part of $Ax-b$. Hoffman's work has been well-recognized and extensively studied  and extended by many authors, including Robinson~\cite{54}, Mangasarian~\cite{42}, Auslender \& Crouzeix \cite{1}, Pang \cite{52}, Lewis \& Pang \cite{41}, Klatte \& Li \cite{34}, Jourani \cite{jourani} , Abassi \& Th\'era \cite{Abassi-Thera1,Abassi-Thera}, Ioffe \cite{ioffe-book} and many others. 
The abundant literature on error bounds shows that this theory has important applications in the sensitivity analysis of linear/integer programs (cf. \cite{Rob73,Rob77}), the convergence analysis of descent methods for linearly constrained minimization (cf. \cite{Gul92,HLu,IuD90,TsL92,TsB93}), the feasibility problem of finding a point in the intersection of finitely many closed convex sets (cf. \cite{5,6,7,BurDeng02}) and the domain of image reconstruction (cf. \cite{16}). The reader is referred to bibliographies \cite{Aze03,Aze06,3,7,CKLT,18,DL,22,ioffe-JAMS-1,ioffe-JAMS-2,Kru15,Luke,44,48,penot-book,55,ZTY-JCA,YZ} and references therein for more details on theoretical research and practical applications.


Due to inaccurate data originated from real-world problems, it is necessary to study the behavior of error bounds when perturbing the data of the problem, and from viewpoints of practical and theoretical interest,   the stability of error bounds under perturbations of the data is expected.   In 1994, Luo and Tseng \cite{LT} studied  the perturbation analysis of a condition number for linear system. In 1998, Deng \cite{D} further discussed the perturbation analysis of a condition number for convex inequality systems. In 2002, Az\'e and Corvellec \cite{2} considered the sensitivity analysis of Hoffman's constants to semi-infinite linear constraint systems. In 2010, Ngai, Kruger and Th\'era \cite{45} studied the stability of error bounds for semi-infinite convex constraint systems in a Euclidean space and established subdifferential characterizations for such stability. Further, Kruger, Ngai \& Th\'era \cite{38} considered the stability of error bounds for semi-infinite convex constraint systems in a Banach space and gave its dual characterizations with respect to perturbations of functions entering in the system. In 2012, by relaxing the convexity assumption, Zheng and Wei \cite{ZW} discussed the stability of error bounds for a semi-infinite constraint system defined by quasi-subsmooth (not necessarily convex) inequalities and provided Clarke subdifferential characterizations of the stability. In 2018, Kruger, Lop\'ez and Th\'era \cite{MP2018} further developed the subdifferential characterizations of the stability of error bounds given in \cite{38,45}. Recently, in terms of directional derivatives of convex functions, the authors \cite{WTY} proved characterizations of  the stability of error bounds for convex inequalities constraint systems in a Euclidean space. In this paper, we continue the study given in \cite{WTY}  on the stability of error bounds for semi-infinite convex constraint systems in a Banach space and 
with an emphasis to primal characterizations of this  stability when perturbing the data of the problem under consideration. It is proved that to verify the stability of error bounds for such  systems,  to some degree, can be equivalent to solving a class of minimax optimization problems determined by directional derivatives of the component functions entering in the system.  When these stability results are  applied to the sensitivity analysis of Hoffman's constants for semi-infinite linear systems, as a result we obtain primal criteria for the Hoffman's constants  to be uniformly bounded under perturbations on the problem data.

The paper is organized as follows. In Section 2, we give some definitions and preliminary results used in the paper. Section 3 is devoted to main results on stability of local and global error bounds for semi-infinite convex constraint systems in a Banach space. In terms of directional derivatives of component functions, we provide primal characterizations of the stability of error bounds under small perturbations (see \cref{th3.1,th3.2}). In Section 4, we apply these stability results to the sensitivity analysis of  the Hoffman's constant  for semi-infinite linear systems, and provide primal criteria for Hoffman's constants to be uniformly bounded under perturbations (see  \cref{th4.1}). Concluding remarks and perspectives of this paper are given in Section 5.

\section{Preliminaries}

Let $\mathbb{X}$  be a Banach space with topological dual $\mathbb{X}^*$. We denote by $\mathbb{B}_{\mathbb{X}}$ the closed unit ball of $\mathbb{X}$ and for  any $x\in \mathbb{X}$ and $\delta>0$,  by $\mathbf{B}(x,\delta)$ the open ball with center $x$ and radius $\delta$.

Let $D$ be a subset of $\mathbb{X}$. We denote by ${\rm int}(D)$ and ${\rm bdry}(D)$ the interior and the boundary of $D$, respectively. The distance from a point $x\in\mathbb{X}$ to the set $D$ is given  by $$\mathbf{d}(x, D):=\inf_{u\in D}\|x-u\|,$$  and we use the convention $\mathbf{d}(x, D) =+\infty$ whenever $D=\emptyset$.

Let $\varphi:\mathbb{X}\rightarrow\mathbb{R}\cup\{+\infty\}$ be a proper lower semicontinuous and convex function and $\bar x\in {\rm dom}(\varphi):=\{x\in \mathbb{X}: \varphi(x)<+\infty\}$. For any $h\in \mathbb{X}$, we denote by $d^+\varphi(\bar x,h)$ the \textit{directional derivative} of $\varphi$ at $\bar x$ along  the direction $h$, which is defined by
\begin{equation}\label{2.1}
d^+\varphi(\bar x,h):=\lim\limits_{t\rightarrow0^+}\frac{\varphi(\bar x+th)-\varphi(\bar x)}{t}.
\end{equation}
According for example to \cite{P}  we know that
$$
t\mapsto \frac{\varphi(\bar x+th)-\varphi(\bar x)}{t}
$$
is nonincreasing as $t\rightarrow 0^+$, and thus
\begin{equation}\label{2.2}
d^+\varphi(\bar x,h)=\inf_{t>0}\frac{\varphi(\bar x+th)-\varphi(\bar x)}{t}.
\end{equation}

We denote by $\partial \varphi(\bar x)$ the \textit{subdifferential} of $\varphi$ at $\bar x$ which is defined by
$$
\partial \varphi(\bar x):=\{x^*\in  \mathbb{X}^*: \langle x^*,x-\bar x\rangle\leq \varphi(x)-\varphi(\bar x)\ \ {\rm for\ all} \ x\in  \mathbb{X}\}.
$$
It is known from \cite[Proposition 2.24]{P} that if $\varphi$ is continuous at $\bar x$ then $\partial \varphi(\bar x)\not=\emptyset$,
\begin{equation}\label{2.3}
 \partial \varphi(\bar x)=\{x^*\in \mathbb{X}^*: \langle x^*,h\rangle\leq d^+\varphi(\bar x,h)\ \  {\rm for \ all} \ h\in  \mathbb{X}\},
\end{equation}
and
\begin{equation}\label{2.4}
d^+\varphi(\bar x,h)=\sup\{ \langle x^*,h\rangle: x^* \in\partial \varphi(\bar x)\}.
\end{equation}

Given a mapping $\phi:\mathbb{X}\rightarrow \mathbb{Y}$ from $\mathbb{X}$ to a \textbf{normed linear space } $\mathbb{Y}$, we denote by
$$
{\rm Lip}(\phi):=\sup_{u,v\in \mathbb{X},u\not=v}\frac{\|\phi(u)-\phi(v)\|_{\mathbb{Y}}}{\|u-v\|_{\mathbb{X}}}.
$$
the Lipschitz constant of $\phi$.

\setcounter{equation}{0}

\section{Stability of Error Bounds for Semi-infinite Convex Constraint Systems}

This section is devoted to the stability of local and global error bounds for semi-infinite convex constraint systems and aims to establish primal characterizations in terms of directional derivatives. We first recall some notations and definitions. 

Let $T$ be a compact, possibly infinite Hausdorff space, and consider the space $ \mathcal{C}(T\times  \mathbb{X}, \mathbb{R})$ of continuous functions on $T\times  \mathbb{X}$.  Let $F\in   \mathcal{C}(T\times  \mathbb{X}, \mathbb{R})$ be such that $F(t,\cdot)$ is convex for each $t\in T$. By a  semi-infinite convex constraint system in $\mathbb{X}$, we mean  the problem of finding $x\in\mathbb{X}$ satisfying:
\begin{equation}\label{4.1}
  f_t(x)\leq 0\ \ {\rm for \ all} \ t\in T,
\end{equation}
where $f_t:= F(t,\cdot)$ for all $t\in T$. Let $S_F$ denote the set of solutions to \eqref{4.1}; that is,
\begin{equation}\label{4.2}
  S_F:=\{x\in  X: f_t(x)\leq 0\  \  {\rm for\ all} \ t\in T\}.
\end{equation}
Set
\begin{equation}\label{4.3}
f(x):=\max\{f_t(x): t\in T\} \ \ {\rm and} \ \ T_f(x):=\{t\in T: f_t(x)=f(x)\},\ \ {\rm for\ each} \ x\in \mathbb{X}.
\end{equation}
Recall that system \eqref{4.1} is said to have a \textit{global error bound}, if there exists a real $c\in(0,+\infty)$ such that
\begin{equation}\label{4.4}
 c \mathbf{d}(x, S_F)\leq [f(x)]_+, \ \ \forall x\in \mathbb{X},
\end{equation}
where $[f(x)]_+:=\max\{f(x), 0\}$. We denote by
$$
{\rm Er} F:=\inf_{x\not\in S_F}\frac{f(x)}{\mathbf{d}(x, S_F)}
$$
the \textit{global error bound modulus} of \eqref{4.1}.

Let $\bar x\in{\rm bdry}(S_F)$. System \eqref{4.1} is said to have a \textit{local error bound} at $\bar x$, if there exist $c,\delta\in(0,+\infty)$ such that
\begin{equation}\label{4.5}
  c \mathbf{d}(x, S_F)\leq [f(x)]_+, \ \ \forall x\in \mathbf{B}(\bar x,\delta).
\end{equation}
We denote by
$$
{\rm Er} F(\bar x):=\liminf_{x\rightarrow \bar x, x\not\in S_F}\frac{f(x)}{\mathbf{d}(x, S_F)}
$$
the \textit{local error bound modulus} of \eqref{4.1} at $\bar x$.

\vskip 2mm

Now, we are in a position to present main work on stability of local and global error bounds for the semi-infinite convex constraint system in a Banach space.  The following theorem gives primal characterizations on the stability of local error bounds for the semi-infinite convex constraint system \eqref{4.1}.

\begin{theorem}\label{th3.1}
Let $F\in \mathcal{C}(T\times  \mathbb{X}, \mathbb{R})$ be such that $f_t:=F(t,\cdot)$ is convex for all $t\in T$,  $\bar x\in  {\rm bdry}(S_F)$ and $f,T_f(x)$ be given as in \eqref{4.3}. Then the following statements are equivalent:
\begin{itemize}
\item[\rm(i)] The optimal value of  minimax problem $
\inf\limits_{\|h\|=1}\sup\limits_{t\in T_f(\bar x)}d^+f_t(\bar x, h)
$ is non-zero.


\item[\rm(ii)] There exist $c,\varepsilon>0$ such that if
\begin{eqnarray*}
&G\in \mathcal{C}(T\times  \mathbb{X}, \mathbb{R}),  g_t(x):=G(t,x), g_t \ {\it is \ convex} ;\\
&g(x):=\sup\limits_{t\in T}g_t(x),  T_g(x):=\{t\in T: g_t(x)=g(x)\};\\
&g(\bar x)=0;\\
&T_g(\bar x)\subseteq T_f(\bar x)\ {\it whenever}\  \inf\limits_{\|h\|=1}\sup\limits_{t\in T_f(\bar x)}d^+f_t(\bar x, h)<0;\\
&T_f(\bar x)\subseteq T_g(\bar x)\ {\it whenever}\  \inf\limits_{\|h\|=1}\sup\limits_{t\in T_f(\bar x)}d^+f_t(\bar x, h)>0;\\
&\limsup\limits_{x\rightarrow\bar x}\displaystyle\frac{|f_t(x)-g_t(x)-(f_t(\bar x)-g_t(\bar x))|}{\|x-\bar x\|}\leq\varepsilon, \ \forall t\in T_f(\bar x)\cap T_g(\bar x),
\end{eqnarray*}
one has ${\rm Er} G(\bar x)\geq c$.
\item[\rm(iii)] There exist $c,\varepsilon>0$ such that for all $u^*\in \mathbb{B}_{\mathbb{X}^*}$, one has ${\rm Er} G_{u^*}(\bar x)\geq c$, where $G_{u^*}\in  \mathcal{C}(T\times  \mathbb{X}, \mathbb{R})$ is defined as follows:
    \begin{equation}\label{4.15}
      G_{u^*}(t,x):=f_t(x)+\varepsilon\langle u^*, x-\bar x\rangle,\ \ {\it for\ all} \ (t,x)\in T\times \mathbb{X}.
    \end{equation}
\end{itemize}
\end{theorem}

\vskip 2mm

The following theorem provides  primal characterizations of the stability of the global error bounds for semi-infinite convex constraint system \eqref{4.1}.

\begin{theorem}\label{th3.2}
Let $F\in \mathcal{C}(T\times  \mathbb{X}, \mathbb{R})$ be such that $f_t:=F(t,\cdot)$ is convex for all $t\in T$ and $f,T_f(x)$ be given as in \eqref{4.3}.
Then the following statements are equivalent:
\begin{itemize}
\item[\rm(i)] There exists $\tau\in (0, +\infty)$ such that
\begin{equation}\label{4.16}
\inf\left\{\Big|\inf_{\|h\|=1}\sup_{t\in T_f(\bar x)}d^+f_t(\bar x,h)\Big|: \bar x\in{\rm bdry}(S_F)\right\}\geq\tau,
\end{equation}
and
\begin{equation}\label{3.20a}
\liminf_{k\rightarrow\infty}\Big|\inf_{\|h\|=1}\sup_{t\in T_f(z_k)}d^+f_t(z_k,h)\Big|\geq\tau
\end{equation}
holds for all $\{(z_k,x_k)\}\subseteq {\rm int}(S_F)\times{\rm bdry}(S_F)$  with $\lim\limits_{k\rightarrow\infty}\frac{f(z_k)-f(x_k)}{\|z_k-x_k\|}=0$.
\item[\rm(ii)] There exist $c,\varepsilon>0$ such that if
\begin{eqnarray*}
&G\in \mathcal{C}(T\times  \mathbb{X}, \mathbb{R}), g_t(x):=G(t,x), g_t \ {\it is\ convex};\\
&g(x):=\sup\limits_{t\in T}g_t(x),  T_g(x):=\{t\in T: g_t(x)=g(x)\};\\                                                                    &\big\{z\in{\rm bdry}(S_F): f_t(z)=g_t(z) \  \ {\it for \ all } \ t\in T\big\}\not=\emptyset;  \\
&\sup\limits_{t\in T}{\rm Lip}(f_t-g_t)<\varepsilon;\\
&T_g(x)\subseteq T_f(x)\ \ {\it whenever} \ \ \inf\limits_{\|h\|=1}\sup\limits_{t\in T_f(x)}d^+f_t(x, h)<0;\\
&T_f(x)\subseteq T_g(x)\ \ {\it whenever} \ \ \inf\limits_{\|h\|=1}\sup\limits_{t\in T_f(x)}d^+f_t(x, h)>0,
\end{eqnarray*}
one has ${\rm Er}G\geq c$.
\item[\rm(iii)] There exist $c,\varepsilon>0$ such that for all $\bar x\in{\rm bdry}(S_F)$ and $u^*\in \mathbb{B}_{\mathbb{X}^*}$, one has ${\rm Er} G_{u^*,\bar x}\geq c$, where $G_{u^*,\bar x}\in  \mathcal{C}(T\times  \mathbb{X}, \mathbb{R})$ is defined as follows:
    \begin{equation}\label{4.23}
      G_{u^*,\bar x}(t,x):=f_t(x)+\varepsilon\langle u^*, x-\bar x\rangle,\ \ {\it for \ all } \ (t, x)\in T\times \mathbb{X}.
    \end{equation}
\end{itemize}
\end{theorem}

Next, our main attention will be drawn to the complete proofs of \cref{th3.1,th3.2}. It is easy to verify that error bounds of system \eqref{4.1} is equivalent to error bounds for $f$ given as in \eqref{4.3}. Thus to prove \cref{th3.1,th3.2}, it is necessary to study the stability of error bounds for the continuous function $f$. Recall that the error bound property for $f$ is defined by the inequality:
\begin{equation}\label{1.1}
  c\mathbf{d}(x, S_f)\leq[f(x)]_+
\end{equation}
where $S_f:=\{x\in \mathbb{X}: f(x)\leq 0\}$ denotes the lower level set of $f$. If inequality \eqref{1.1} holds for some $c>0$ and all $x\in \mathbb{X}$, then $f$ is said to have a global error bound. In this case, we denote by
\begin{equation}\label{1.2}
  {\rm Er}f:=\inf_{f(x)>0}\frac{f(x)}{\mathbf{d}(x,S_f)}
\end{equation}
the global error bound modulus of $f$. Given $\bar x\in X$ with $f(\bar x)=0$, if inequality \eqref{1.1} holds for some $c>0$  at all $x$ near to $\bar x$, then $f$ is said to have a local error bound at $\bar x$. In this case, we denote
\begin{equation}\label{1.3}
  {\rm Er}f(\bar x):=\liminf_{x\rightarrow \bar x, f(x)>0}\frac{f(x)}{\mathbf{d}(x,S_f)}
\end{equation}
the local error bound modulus of $f$ at $\bar x$.

\subsection{Stability of Error Bounds for $f$}

We first prove the following proposition on primal characterizations of error bounds for the continuous convex function $f$ given as in \eqref{4.3}. This is a key tool to prove main results in this paper.

\begin{proposition}\label{pro3.1}
Let $f$ be given in \eqref{4.3} such that the lower level set $S_f$ is nonempty.
\begin{itemize}
\item[\rm (i)] $f$ has a global error bound if and only if
\begin{equation*}\label{3.6a}
  \beta(f):=\inf_{f(x)>0}\left(-\inf_{\|h\|=1}d^+ f(x,h)\right)>0.
\end{equation*}
Moreover, ${\rm Er} f=\beta(f)$.
\item[\rm (ii)] $f$ has a local error bound at $\bar x\in{\rm bdry}(S_{f})$ if and only if
\begin{equation*}\label{3.6}
  \beta(f,\bar x):=\liminf_{x\rightarrow\bar x, f(x)>0}\left(-\inf_{\|h\|=1}d^+ f(x,h)\right)>0.
\end{equation*}
Moreover, ${\rm Er} f(\bar x)=\beta(f,\bar x)$.
\end{itemize}
\end{proposition}

{\bf Proof.} We first prove that
\begin{equation}\label{3.8}
  -\inf_{\|h\|=1}d^+ f(x,h)=\mathbf{d}(0,\partial f(x)),\ \ \forall x: f(x)>0.
\end{equation}

Let $x\in \mathbb{X}$ be such that $f(x)>0$. Then $\partial f(x)\not=\emptyset$ and $0\not\in\partial f(x)$ (thanks to $f(x)>0$ and $S_{f}\not=\emptyset$). We denote $r:=\mathbf{d}(0,\partial f(x))>0$.
Take any $\varepsilon>0$ sufficiently small  in order to have \\ $(r-\varepsilon) \mathbb{B}_{\mathbb{X}^*}\cap\partial f(x)=\emptyset$.  Noting that $\mathbb{B}_{\mathbb{X}^*}$ and $\partial f(x)$ are weak$^*$-compact, it follows from the  separation theorem in \cite[Theorem 3.4]{Rudin} that there exist $h_0\in \mathbb{X}$ with $\|h_0\|=1$ and $\gamma_1\in\mathbb{R}$ such that
$$
\sup\{\langle x^*, h_0\rangle: x^*\in\partial f(x)\}<\gamma_1<\inf\{\langle x^*, h_0\rangle: x^*\in (r-\varepsilon)\mathbb{B}_{\mathbb{X}^*})\}=-(r-\varepsilon).
$$
This and \eqref{2.4} imply that
\begin{equation}\label{3.9a}
 \inf_{\|h\|=1}d^+f(x ,h)\leq d^+f(x ,h_0)<-(r-\varepsilon)\rightarrow -r\ \ (\rm as \ \varepsilon\rightarrow 0^+).
\end{equation}
On the other hand, let $\varepsilon>0$. We can select $u^*_{\varepsilon}\in (r+\varepsilon) \mathbb{B}_{\mathbb{X}^*}\cap\partial f(x)$. Then for any $h\in \mathbb{X}, \|h\|=1$, one has
$$
d^+f(x ,h)\geq \langle u^*_{\varepsilon}, h\rangle\geq -(r+\varepsilon),
$$
and consequently
$$
\inf_{\|h\|=1}d^+f(x ,h)\geq -(r+\varepsilon).
$$
Letting $\varepsilon\rightarrow 0^+$, one has \eqref{3.8} holds by \eqref{3.9a}.

By virtue of \cite[Theorem 1 and Theorem 22]{38}, the following equalities on error bound moduli for $f$ hold:
\begin{equation}\label{3-8a}
  {\rm Er} f=\inf_{f(x)>0}\mathbf{d}(0,\partial f(x))\ \ {\rm and} \ \  {\rm Er} f(\bar x)=\liminf_{x\rightarrow\bar x, f(x)>0}\mathbf{d}(0,\partial f(x)).
\end{equation}
This and \eqref{3.8} imply that the conclusions in (i) and (ii) follow. The proof is complete.\qed

\vskip 2mm

The following proposition is to present primal criteria of the local error bound for the function $f$.

\begin{proposition}\label{pro3.2}
Let $f$ be given as in \eqref{4.3}, $\bar x\in {\rm bdry}(S_{F})$ and $\inf_{\|h\|=1}d^+f(\bar x,h)\not=0$.
\begin{itemize}
\item[\rm(i)] If $\inf_{\|h\|=1}d^+f(\bar x, h)>0$, then $S_{F}=\{\bar x\}$ and ${\rm Er} f(\bar x)\geq \inf_{\|h\|=1}d^+f(\bar x,h)$.
\item[\rm(ii)] If $\inf_{\|h\|=1}d^+f(\bar x, h)<0$, then ${\rm Er} f(\bar x)\geq -\inf_{\|h\|=1}d^+f(\bar x,h)$.
\end{itemize}
\end{proposition}

{\bf Proof.} (i) It is easy to verify that $S_F=S_f$ and $f(\bar x)=0$. Let $x\in \mathbb{X}\backslash\{\bar x\}$.  By \eqref{2.2}, one has
\begin{eqnarray*}
f(x)-f(\bar x)&=&f\Big(\bar x+\|x-\bar x\|\frac{x-\bar x}{\|x-\bar x\|}\Big)-f(\bar x)\\
&\geq& d^+f\Big(\bar x,\frac{x-\bar x}{\|x-\bar x\|}\Big)\|x-\bar x\|>0.
\end{eqnarray*}
This implies that $S_{f}=\{\bar x\}$ and ${\rm Er} f(\bar x)\geq \inf_{\|h\|=1}d^+f(\bar x,h)$.

(ii) By virtue of \eqref{3.8} and \eqref{3-8a}, one has
$$
{\rm Er} f(\bar x)=\liminf_{x\rightarrow\bar x, f(x)>0}\mathbf{d}(0,\partial f(x))\geq \mathbf{d}(0,\partial f(\bar x))=-\inf_{\|h\|=1}d^+f(\bar x,h).
$$
Thus the conclusion holds. The proof is complete.\hfill$\Box$

\begin{remark} The inequalities in \cref{pro3.2} may hold strictly. For example, let $\mathbb{X}:=\mathbb{R}$ and $F(t,x)\equiv 0$ for all $(t, x)\in T\times\mathbb{R}$. Then $f$ has the local error bound at any $\bar x\in \mathbb{R}$ and ${\rm Er} f(\bar x)=+\infty$. However, one can verify that $\inf_{\|h\|=1}d^+f(\bar x, h)=0$, and thus ${\rm Er} f(\bar x)> \inf_{\|h\|=1}d^+f(\bar x,h)$.
\end{remark}

The following theorem shows that the condition $$\inf_{\|h\|=1}d^+f(\bar x,h)\not=0$$ is sufficient and necessary for the stability of local error bounds for $f$ at $\bar x$.

\begin{theorem}\label{th3.3}
Let $f$ be  as \eqref{4.3} and $\bar x\in {\rm bdry}(S_{F})$. Then the following statements are equivalent:
\begin{itemize}
\item[\rm (i)] $\inf_{\|h\|=1}d^+f(\bar x, h)\not=0$.
\item[\rm (ii)] There exist $c,\varepsilon>0$ such that for all $g\in  \mathcal{C}(\mathbb{X},\mathbb{R})$ satisfying $\bar x\in S_g$ and
\begin{equation}\label{3.9}
\limsup_{x\rightarrow \bar x}\frac{|(f(x)-g(x))-(f(\bar x)-g(\bar x))|}{\|x-\bar x\|}\leq \varepsilon,
\end{equation}
one has ${\rm Er} g(\bar x)\geq c$.
\item[\rm (iii)] There exist $c,\varepsilon>0$ such that for all $u^*\in \mathbb{B}_{\mathbb{X}^*}$, one has ${\rm Er} g_{u^*, \varepsilon}(\bar x)\geq c$, where $g_{u^*, \varepsilon}(\cdot):=f(\cdot)+\varepsilon\langle u^*, \cdot-\bar x\rangle$.
\end{itemize}
\end{theorem}

{\bf Proof.} We denote
$$
\gamma:=\inf_{\|h\|=1}d^+f(\bar x, h).
$$

(i) $\Rightarrow$ (ii): Let $\varepsilon>0$ be such that $\varepsilon<|\gamma|$ and $c:=|\gamma|-\varepsilon$. Take any $g\in  \mathcal{C}(\mathbb{X},\mathbb{R})$ satisfying $\bar x\in S_g$ and \eqref{3.9}. If $\gamma>0$, then for any $h\in \mathbb{X}$ with $\|h\|=1$, one has
$$
d^+g(\bar x, h)\geq d^+f(\bar x, h)-\varepsilon,
$$
and consequently
$$
\inf_{\|h\|=1}d^+g(\bar x, h)\geq\inf_{\|h\|=1}d^+f(\bar x, h)-\varepsilon= \gamma-\varepsilon.
$$
Using \cref{pro3.2}, one has ${\rm Er} g(\bar x)\geq \gamma-\varepsilon$.

If $\gamma<0$, then for any $h\in \mathbb{X}$, one has
$$
d^+g(\bar x, h)\leq d^+f(\bar x, h)+\varepsilon,
$$
and thus
$$
\inf_{\|h\|=1}d^+g(\bar x, h)\leq\inf_{\|h\|=1}d^+f(\bar x, h)+\varepsilon=\gamma+\varepsilon.
$$
Using \cref{pro3.2} again, one can get ${\rm Er} g(\bar x)\geq -\gamma-\varepsilon=c$. Hence (ii) holds.

Note that (ii) $\Rightarrow$ (iii) follows immediately and we next prove (iii) $\Rightarrow$ (i).

Let $\varepsilon>0$. Suppose on the contrary that $\gamma=0$. Then there is a sequence $\{h_k\}\subseteq \mathbb{X}$ such that $\|h_k\|=1$ and
$$
\alpha_k:=d^+f(\bar x,h_k)\rightarrow \gamma=0\ \ {\rm as} \ k\rightarrow\infty.
$$
Without loss of generality, we can assume that $|\alpha_k|<\varepsilon$ holds for all $k$ (considering sufficiently large $k$ if necessary). Then for any $k$, by the Hahn-Banach theorem, there is $h^*_k\in \mathbb{X}^*$ such that
$$
\|h^*_k\|=1\ \ {\rm and} \ \ \langle h_k^*, h_k\rangle=\|h_k\|.
$$
We now consider the following function:
$$
g_{\varepsilon}(x):=f(x)+\varepsilon \langle h_k^*, x-\bar x\rangle, \forall x\in \mathbb{X}.
$$
Note that $\gamma=0$ and thus for any $x\not=\bar x$, one has $f(x)\geq f(\bar x)$. By the definition of  a directional derivative, there exists $\delta_k\rightarrow 0^+$ such that
\begin{equation}\label{3.10}
  f(\bar x+\delta_kh_k)<f(\bar x)+(\varepsilon+\alpha_k)\delta_k=\inf_{x\in \mathbb{X}}f(x)+(\varepsilon+\alpha_k)\delta_k.
\end{equation}
By virtue of  Ekeland's variational principle(cf. \cite{E1974}), there is $z_k\in \mathbb{X}$ such that
$$
\|z_k-(\bar x+\delta_kh_k)\|<\frac{\delta_k}{2}, f(z_k)\leq f(\bar x+\delta_kh_k)
$$
and
\begin{equation}\label{3.11}
 f(x)+2(\varepsilon+\alpha_k)\|x-z_k\|>f(z_k),\ \forall x: x\not=z_k.
\end{equation}
This implies that $z_k\rightarrow \bar x$, $g_{\varepsilon}(\bar x)=f(\bar x)=0$ and
\begin{eqnarray*}
g_{\varepsilon}(z_k)&=&f(z_k)+\varepsilon\langle h_k^*, z_k-\bar x \rangle\\
&\geq& f(\bar x)+\varepsilon\langle h_k^*, z_k-\bar x \rangle\\
&=&\varepsilon\langle h_k^*, z_k-\bar x-\delta_kh_k \rangle+\varepsilon\delta_k\\
&>&\varepsilon\delta_k-\frac{1}{2}\varepsilon\delta_k=\frac{1}{2}\varepsilon\delta_k>0.
\end{eqnarray*}
We claim that
\begin{equation}\label{3.12}
  \inf_{\|h\|=1}d^+g_{\varepsilon}(z_k,h)< 0,
\end{equation}
since otherwise, $\inf_{\|h\|=1}d^+g_{\varepsilon}(z_k,h)\geq 0$, and thus $g_{\varepsilon}(z_k)=\inf_{x\in \mathbb{X}}g_{\varepsilon}(x)$, which contradicts $g_{\varepsilon}(\bar x)=0$.

For any $h\in \mathbb{X}$, $\|h\|=1$ and any $t>0$, by \eqref{3.11}, one can get
\begin{eqnarray*}
\frac{g_{\varepsilon}(z_k+th)-g_{\varepsilon}(z_k)}{t}=\frac{f(z_k+th)-f(z_k)}{t}+\varepsilon\langle h_k^*, h\rangle\geq-2(\varepsilon+\alpha_k)\|h\|-\varepsilon=-5\varepsilon.
\end{eqnarray*}
Thus,
$$
\inf_{\|h\|=1}d^+g_{\varepsilon}(z_k,h)\geq-5\varepsilon.
$$
Combining \eqref{3.12} with \cref{pro3.1} yields   ${\rm Er} g_{\varepsilon}(\bar x)\leq5\varepsilon\rightarrow 0^+$, which contradicts (iii), as $\varepsilon$ is arbitrary. The proof is complete. \hfill$\Box$

\begin{remark}
By \cite[Definition 5]{38}, condition \eqref{3.9} is called as an
\textit{$\varepsilon$-perturbation } of $f$ near $\bar x$. It is shown in \cref{th3.3} that condition ~$\inf_{\|h\|=1}d^+f(\bar x, h)\not=0$ is sufficient and necessary for a local error bound moduli  of all $\varepsilon$-perturbations of $f$ near $\bar x$ to be uniformly bounded. Further, such stability under perturbations on $f$ only requires all the linear $\varepsilon$-perturbations of $f$ near $\bar x$.
\end{remark}

The following corollaries are immediate from \cref{th3.3}.

\begin{corollary}
Let $f$ be given as  in \eqref{4.3} and $\bar x\in {\rm bdry}(S_{F})$. Then
 $\inf\limits_{\|h\|=1}d^+f(\bar x,h)\not=0$ holds if and only if there is an $\varepsilon>0$ such that
$$
\inf\big\{{\rm Er} g(\bar x): g\in  \mathcal{C}(\mathbb{X},\mathbb{R})\ {\it satisfying} \ \bar x\in S_g\ {\it and} \ \eqref{3.9} \big\}>0.
$$
\end{corollary}

\begin{corollary}
Let $f$ be given as in  \eqref{4.3}, $\bar x\in {\rm bdry}(S_{F})$ and $\varepsilon>0$.
\begin{itemize}
\item[\rm(i)] The following inequality holds:
$$
\inf\big\{{\rm Er} g(\bar x): g\in   \mathcal{C}(\mathbb{X},\mathbb{R})\ {\it satisfying} \ \bar x\in S_g\ {\it and} \ \eqref{3.9} \  \big\}\ \geq\left|\inf_{\|h\|=1}d^+f(\bar x,h)\right|-\varepsilon.
$$
\item[\rm(ii)] If $\inf_{\|h\|=1}d^+f(\bar x,h)=0$, then
$$
\inf\big\{{\rm Er} g(\bar x): g\in   \mathcal{C}(\mathbb{X},\mathbb{R})\ {\it satisfying} \ \bar x\in S_g\ {\it and} \ \eqref{3.9} \ \big\}=0.
$$
\end{itemize}
\end{corollary}

The following proposition is to give primal criteria for the stability of the global error bound for the continuous convex function $f$.


\begin{theorem}\label{th3.4}
Let $f$ be given as   in \eqref{4.3}.
\begin{itemize}
\item[\rm (i)] If there is $\tau\in (0, +\infty)$ such that
\begin{equation}\label{3-9}
\inf\left\{\Big|\inf_{\|h\|=1}d^+f(\bar x,h)\Big|: \bar x\in{\rm bdry}(S_f)\right\}\geq\tau,
\end{equation}
then there exist $c,\varepsilon\in (0, +\infty)$ such that for all $g\in \mathcal{C}(T\times  \mathbb{X}, \mathbb{R})$ satisfying
    \begin{equation}\label{3-10}
      S_f\subseteq S_g\  \ {\it and} \ \ {\rm Lip}(f-g)<\varepsilon,
    \end{equation}
one has ${\rm Er} g\geq c$.
\item[\rm (ii)] If there are $c,\varepsilon\in (0, +\infty)$ such that ${\rm Er} g\geq c$ holds for all $g\in  \mathcal{C}(\mathbb{X}, \mathbb{R})$ satisfying
    \begin{equation}\label{3-11}
      {\rm bdry}(S_f)\cap g^{-1}(0)\not=\emptyset\ \  {\it and}  \  \ {\rm Lip}(f-g)<\varepsilon,
    \end{equation}
then there exists $\tau\in (0, +\infty)$ such that \eqref{3-9} holds
\end{itemize}
\end{theorem}

{\bf Proof.} Note that $f$ is continuous and thus ${\rm bdry}(S_f)\subseteq f^{-1}(0)$.

(i)  Suppose that there is some $\bar x\in{\rm bdry}(S_f)$ such that $\inf_{\|h\|=1}d^+f(\bar x, h)>0$. By \cref{pro3.2}, the lower level set $S_f$ is a  singleton and thus $S_f=\{\bar x\}$. This means that the conclusion in (i) follows by using \cref{th3.3} as well as its proof.

Next, we assume that $\inf_{\|h\|=1}d^+f(\bar x, h)\leq 0$ for all $\bar x\in{\rm bdry}(S_f)$. We claim that
\begin{equation}\label{3-12}
  \tau\mathbf{d}(x, S_f)\leq [f(x)]_+,\ \ \forall x\in \mathbb{X}.
\end{equation}
Let $x\in\mathbb{X}\backslash S_f$ and $\nu>0$ sufficiently small. Take $x_{\nu}\in {\rm bdry}(S_{f})$ such that
\begin{equation}\label{3.017a}
  \|x-x_{\nu}\|<(1+\nu)\mathbf{d}(x, S_f).
\end{equation}
By virtue of \cref{pro3.1} and \eqref{3.8}, one has
$$
{\rm Er}f(x_{\nu})=\liminf_{u\rightarrow x_{\nu}, u\not\in S_f}\Big(-\inf_{\|h\|=1}d^+f(u, h)\Big)\geq-\inf_{\|h\|=1}d^+f(x_{\nu}, h)\geq\tau.
$$
Then there exists $\lambda\in(0, 1)$ such that
\begin{eqnarray*}
\tau\leq\frac{f(\lambda x_{\nu}+(1-\lambda)x)}{\mathbf{d}(\lambda x_{\nu}+(1-\lambda)x, S_f)}+\nu\leq\frac{\lambda f(x_{\nu})+(1-\lambda)f(x)}{\mathbf{d}(x, S_f)-\|\lambda x_{\nu}+(1-\lambda)x-x\|}+\nu,
\end{eqnarray*}
and it follows from \eqref{3.017a} that
\begin{eqnarray*}
\tau\leq\frac{(1-\lambda)f(x)}{(1-\lambda(1+\nu))\mathbf{d}(x, S_f)}+\nu=\frac{f(x)}{\mathbf{d}(x, S_f)}+\frac{\nu \lambda f(x)}{\mathbf{d}(x, S_f)}+\nu\rightarrow \frac{f(x)}{\mathbf{d}(x, S_f)}\ {\rm as\ }\nu\rightarrow 0^+.
\end{eqnarray*}
This means that \eqref{3-12} holds.

Let $\varepsilon\in (0,\tau)$ and $g\in   \mathcal{C}(\mathbb{X}, \mathbb{R})$ satisfying \eqref{3-10}. Take any $x\in \mathbb{X}$ such that $g(x)>0$. Then $f(x)>0$ (thanks to $S_f\subseteq S_g$).

We claim that
\begin{equation}\label{3.13}
  \inf_{\|h\|=1}d^+f(x, h)\leq -\tau.
\end{equation}
Granting this, it follows from ${\rm Lip}(f-g)<\varepsilon$ in \eqref{3-10} that
$$
\inf_{\|h\|=1}d^+g(x,h)\leq\inf_{\|h\|=1}d^+f(x,h)+\varepsilon\leq -(\tau-\varepsilon).
$$
By \cref{pro3.2}, one has that ${\rm Er} g\geq \tau-\varepsilon$.

It remains to prove \eqref{3.13}. For any  $n\in\mathbb{N}$, we can take $z_n\in {\rm bdry}(S_f)$ such that
$$
\|x-z_n\|<(1+\frac{1}{n})\mathbf{d}(x,S_f).
$$
Then \eqref{3-12} implies that
$$
f(x)\geq \tau \mathbf{d}(x, S_f)>\frac{n}{n+1}\tau \|x-z_n\|.
$$
For any $t\in (0,1)$, one has
$$
f(x+t(z_n-x))\leq tf(z_n)+(1-t)f(x),
$$
and thus
$$
\frac{f(x+t(z_n-x))-f(x)}{t}\leq -f(x)\leq-\frac{n}{n+1}\tau\|x-z_n\|.
$$
This implies that
$$
\inf_{\|h\|=1}d^+f(x, h)\leq d^+f\Big(x,\frac{z_n-x}{\|x-z_n\|}\Big)\leq-\frac{n}{n+1}\tau.
$$
Hence
$$
\inf_{\|h\|=1}d^+f(x, h)\leq-\frac{n}{n+1}\tau.
$$
Letting $n\rightarrow\infty$, one has that \eqref{3.13} holds.

(ii) Suppose on  the contrary that there is a sequence $\{x_k\}\subseteq {\rm bdry}(S_f)$ such that
$$
\alpha_k:=\inf_{\|h\|=1}d^+f(x_k, h)\rightarrow 0^- \ \ {\rm as} \ k\rightarrow \infty.
$$
Let $\varepsilon>0$ and take sufficiently large $k$ such that
\begin{equation}\label{3.14}
  \frac{3}{2}\alpha_k+\frac{\varepsilon}{2}>0.
\end{equation}
Note that for any $x\not=x_k$, one has
$$
\frac{f(x)-f(x_k)}{\|x-x_k\|}=\frac{f\big(x_k+\|x-x_k\|\cdot\frac{x-x_k}{\|x-x_k\|}\big)-f(x_k)}{\|x-x_k\|}\geq d^+f\Big(x_k,\frac{x-x_k}{\|x-x_k\|}\Big)\geq\alpha_k,
$$
and consequently
\begin{equation}\label{3.15}
  f(x)-\alpha_k\|x-x_k\|\geq f(x_k),\ \ \forall x\in \mathbb{X}.
\end{equation}
By the definition of $\alpha_k$, there exists $h_k\in \mathbb{X}$ such that
\begin{equation}\label{3.16}
\|h_k\|=1 \ \ {\rm and} \ \  d^+f(x_k,h_k)<\alpha_k+\frac{\varepsilon}{2}.
\end{equation}
Then one can choose $r_k\rightarrow 0^+$ (as $k\rightarrow \infty$) such that
\begin{equation}\label{3.17}
  f(x_k+r_kh_k)<f(x_k)+(\alpha_k+\varepsilon)r_k.
\end{equation}
This and \eqref{3.15} imply that
$$
  f(x_k+r_kh_k)-\alpha_k\|x_k+r_kh_k-x_k\|<\inf_{x\in X}(f(x)-\alpha_k\|x-x_k\|)+\varepsilon r_k.
$$
Applying  Ekeland's variational principle (cf \cite{E1974}), there exists $y_k\in \mathbb{X}$ such that
\begin{equation}\label{3.18}
  \|y_k-(x_k+r_kh_k)\|<\frac{r_k}{2}, f(y_k)-\alpha_k\|y_k-x_k\|\leq  f(x_k+r_kh_k)-\alpha_k r_k,
\end{equation}
and
\begin{equation}\label{3.19}
  f(x)-\alpha_k\|x-x_k\|+2\varepsilon\|x-y_k\|>f(y_k)-\alpha_k\|y_k-x_k\|,\ \ \forall x: x\not=y_k.
\end{equation}
Then
$$
\|y_k-x_k\|> r_k-\frac{r_k}{2}=\frac{r_k}{2}\ \ {\rm and}\ \ \|y_k-x_k\|< r_k+\frac{r_k}{2}=\frac{3}{2}r_k,
$$
and thus $y_k\not=x_k$.

For any $h_k$, by  the Hahn-Banach theorem, there exists $h^*_k\in \mathbb{X}^*$ such that
$$
\|h^*_k\|=1\ \ {\rm and} \ \ \langle h_k^*, h_k\rangle=\|h_k\|.
$$
We consider the following function $g_{\varepsilon}$:
$$
g_{\varepsilon}(x):=f(x)+\varepsilon\langle h_k^*, x-x_k\rangle\ \ \forall x\in \mathbb{X}.
$$
Then $g_{\varepsilon}\in  \mathcal{C}(\mathbb{X}, \mathbb{R})$. By virtue of \eqref{3.14}, \eqref{3.15}, \eqref{3.16} and \eqref{3.19}, one has
\begin{eqnarray*}
g_{\varepsilon}(y_k)=f(y_k)+\varepsilon\langle h_k^*, y_k-x_k\rangle&=&f(y_k)+\varepsilon\langle h_k^*, y_k-(x_k+r_kh_k)\rangle+\varepsilon r_k\\
&\geq&\alpha_k\|y_k-x_k\|-\varepsilon\|y_k-(x_k+r_kh_k)\|+\varepsilon r_k\\
&\geq&\alpha_k\cdot\frac{3}{2}r_k+\frac{\varepsilon}{2} r_k>0.
\end{eqnarray*}
If $\inf_{\|h\|=1}d^+g_{\varepsilon}(y_k, h)\geq 0$, then for any $x\not=y_k$, one has
$$
g_{\varepsilon}(x)-g_{\varepsilon}(y_k)\geq d^+g_{\varepsilon}\big(y_k,\frac{x-y_k}{\|x-y_k\|}\big)\|x-y_k\|\geq \inf_{\|h\|=1}d^+g_{\varepsilon}(y_k, h)\|x-y_k\|\geq 0.
$$
This together with $g_{\varepsilon}(y_k)>0$ implies that  $S_{g_{\varepsilon}}=\emptyset$ and thus ${\rm Er} g_{\varepsilon}=0$, which contradicts (ii).

Next, we consider $\inf_{\|h\|=1}d^+g_{\varepsilon}(y_k, h)<0$. For any $h\in \mathbb{X}$ with $\|h\|=1$ and any  $t>0$, by \eqref{3.19}, one has
\begin{eqnarray*}
\frac{g_{\varepsilon}(y_k+th)-g_{\varepsilon}(y_k)}{t}&=&\frac{f(y_k+th)-f(y_k)}{t}+\varepsilon\langle h_k^*, h\rangle\\
&\geq&\frac{1}{t}\big(\alpha_k\|y_k+th-x_k\|-\alpha_k\|y_k-x_k\|-2\varepsilon \|y_k+th-y_k\|\big)+\varepsilon\langle h_k^*, h\rangle\\
&\geq&\alpha_k-2\varepsilon-\varepsilon,
\end{eqnarray*}
and consequently
$$
0>\inf_{\|h\|=1}d^+g_{\varepsilon}(y_k, h)\geq\alpha_k-2\varepsilon-\varepsilon\geq-4\varepsilon.
$$
Applying \cref{pro3.1}, one derives that ${\rm Er} (g_{\varepsilon})\leq 4\varepsilon$, which contradicts (ii), as $\varepsilon$ is arbitrary. The proof is complete. \hfill$\Box$

The following corollary follows from the proof of \cref{th3.4}.
\vskip 2mm

\begin{corollary}
Let $f$ be given  as in \eqref{4.3}. Then
\begin{equation}
{\rm Er} f\geq \inf\left\{\Big|\inf_{\|h\|=1}d^+f(\bar x,h)\Big|: \bar x\in{\rm bdry}(S_f)\right\}.
\end{equation}
\end{corollary}

\begin{remark}\label{re3-2}
\cref{th3.4} shows that \eqref{3-9} is a sufficient condition to ensure the stability of the global error bound for $f$  as said in (i) and meanwhile it is a necessary condition  for the stability of the global error bound for $f$  as said in (ii). However, condition \eqref{3-9} may not be sufficient for the stability of the global error bound for $f$  as in (ii). For example, let $\mathbb{X}:=\mathbb{R}$ and $F(t,x):=e^x-1$ for all $(t,x)\in T\times\mathbb{R} $. Then $f(x)=e^x-1$, $S_f=(-\infty, 0]$, ${\rm bdry}(S_f)=\{0\}$ and $|\inf_{|h|=1}d^+f(0, h)|=1>0$. Let $\varepsilon\in (0,+\infty)$, and consider $g_{\varepsilon}(x):=f(x)-\varepsilon x, \forall x\in\mathbb{R}$. Then one can verify that $g_{\varepsilon}$ has two different zero points $0$ and $\bar x<0$. Further, one can verify that $S_{g_{\varepsilon}}=[\bar x, 0]$ and thus $S_f\not\subseteq S_{g_{\varepsilon}}$. However, for any $x<\bar x$, one has
$$
\frac{g_{\varepsilon}(x)}{\mathbf{d}(x, S_{g_{\varepsilon}})}=\frac{e^x-1-\varepsilon x}{\bar x-x}\rightarrow \varepsilon \ \ {\rm as} \ x\rightarrow -\infty.
$$
This implies that
$$
{\rm Er} g_{\varepsilon}\leq 2\varepsilon,
$$
and consequently the stability of global error bound for $f$  as in (ii) does not hold (as $\varepsilon>0$ is arbitrary).
\end{remark}

The following theorem shows that \eqref{3-9} plus a mild qualification can characterize the stability of the global error bound for $f$  as said in (ii) of \cref{th3.4}.

\begin{theorem}\label{th3.5}
Let $f$ be   as \eqref{4.3}. Then the following statements are equivalent:
\begin{itemize}
\item[\rm (i)] There exists $\tau\in (0, +\infty)$ such that \eqref{3-9} holds and
\begin{equation}\label{3.20}
\liminf_{k\rightarrow\infty}\Big|\inf_{\|h\|=1}d^+f(z_k,h)\Big|\geq\tau, \ \forall \{(z_k,x_k)\}\subseteq {\rm int}(S_f)\times{\rm bdry}(S_f)\ {\it with} \ \lim\limits_{k\rightarrow\infty}\frac{f(z_k)-f(x_k)}{\|z_k-x_k\|}=0.
\end{equation}
\item[\rm (ii)] There exist $c,\varepsilon\in (0, +\infty)$ such that for all $g\in  \mathcal{C}(T\times  \mathbb{X}, \mathbb{R})$ satisfying \eqref{3-11}, one has ${\rm Er} g\geq c$;
\item[\rm (iii)] There exist $c,\varepsilon>0$ such that for any $\bar x\in{\rm bdry}(S_f)$ and any $u^*\in \mathbb{B}_{\mathbb{X}^*}$, one has ${\rm Er} g_{u^*, \varepsilon}\geq c$, where $g_{u^*, \varepsilon}(\cdot):=f(\cdot)+\varepsilon\langle u^*, \cdot-\bar x\rangle$.
\end{itemize}
\end{theorem}

{\bf Proof.} (i)$\Rightarrow$(ii): If there is $\bar x\in{\rm bdry}(S_f)$ such that $\inf_{\|h\|=1}d^+f(\bar x, h)>0$, the conclusion follows by \cref{pro3.2} as well as the proof of (i)$\Rightarrow$ (ii) in \cref{th3.3}.

Next, we consider $\inf_{\|h\|=1}d^+f(\bar x, h)\leq 0$ for all $\bar x\in{\rm bdry}(S_f)$. We first prove the following claim:

{\it There exists $\varepsilon_0>0$ such that for any $x_0\in{\rm bdry}(S_f)$, one has}
\begin{equation}\label{3.21}
\inf\left\{\Big|\inf_{\|h\|=1}d^+f(z_0,h)\Big|: z_0\in \mathbb{X}, f(z_0)\geq-\varepsilon_0\|z_0-x_0\|\right\}\geq\tau.
\end{equation}

Suppose on the contrary that there exist $\varepsilon_k\rightarrow 0^+$, $x_k\in{\rm bdry}(S_f)$ and $z_k\in\mathbb{X}$ such that
\begin{equation}\label{3.22}
  f(z_k)\geq-\varepsilon_k\|z_k-x_k\|\ \ {\rm and} \ \ \Big|\inf_{\|h\|=1}d^+f(z_k,h)\Big|<\tau,\ \ \forall  k.
\end{equation}
Then $f(z_k)\leq 0$  for all $k$, since otherwise, using the proof of \eqref{3.13}, one get $$\big|\inf_{\|h\|=1}d^+f(z_k,h)\big|>\tau,
 $$
which contradicts \eqref{3.22}.

By virtue of \eqref{3-9}, one has $z_k\in S_f\backslash {\rm bdry}(S_f)$ and then \eqref{3.22} implies
$$
0\geq\frac{f(z_k)-f(x_k)}{\|z_k-x_k\|}=\frac{f(z_k)}{\|z_k-x_k\|}\geq-\varepsilon_k.
$$
By \eqref{3.20}, one gets
$$
\liminf_{k\rightarrow\infty}\Big|\inf_{\|h\|=1}d^+f(z_k,h)\Big|>\tau,
$$
which contradicts \eqref{3.22}. Hence the claim holds.

Let $\varepsilon>0$ be such that $\varepsilon<\min\{\varepsilon_0, \tau\}$ and $g\in  \mathcal{C}(\mathbb{X}, \mathbb{R})$ be such that \eqref{3-11} holds. Take any $\bar x\in{\rm bdry}(S_f)\cap g^{-1}(0)$. Then for any $x\in \mathbb{X}$  with $g(x)>0$, one has
$$
f(x)\geq g(x)+(f(\bar x)-g(\bar x))-\varepsilon\|x-\bar x\|>-\varepsilon\|x-\bar x\|.
$$
By \eqref{3.21}, one has
$$
\inf_{\|h\|=1}d^+f(x,h)<-\tau,
$$
and consequently
$$
\inf_{\|h\|=1}d^+g(x,h)<\inf_{\|h\|=1}d^+f(x,h)+\varepsilon<-(\tau-\varepsilon).
$$
Then \cref{pro3.2} implies that ${\rm Er} g\geq \tau-\varepsilon$.

Note (ii)$\Rightarrow$(iii) follows immediately and it remains to prove (iii)$\Rightarrow$(i).

Suppose on the contrary that the conclusions in (i) does not hold. Based on \cref{th3.4} and its proof, it suffices to assume that \eqref{3.20} does not hold for any $\tau>0$. Then, there exists $(z_k,x_k)\in {\rm int}(S_f)\times {\rm bdry}(S_f)$ such that
\begin{equation}\label{3.23}
  \lim_{k\rightarrow\infty}\frac{f(z_k)-f(x_k)}{\|z_k-x_k\|}=0\ \ {\rm and} \ \ \alpha_k:=\inf_{\|h\|=1}d^+f(z_k,h)\rightarrow 0^- \ \ {\rm as}\ k\rightarrow\infty.
\end{equation}
Let $\varepsilon>0$. For any $k$, by the Hahn-Banach theorem,  there is $h_k^*\in \mathbb{X}^*$ such that
$$
\|h_k^*\|=1\ \ {\rm and} \ \ \langle h_k^*, z_k-x_k\rangle=\|z_k-x_k\|.
$$
Then when $k$ is sufficiently large, one has
\begin{equation}\label{3.24}
   \alpha_k+\varepsilon>0\ \  {\rm and} \ \ \frac{f(z_k)-f(x_k)}{\|z_k-x_k\|}+\varepsilon\Big\langle h_k^*, \frac{z_k-x_k}{\|z_k-x_k\|}\Big\rangle>0.
\end{equation}
Consider the function $g_{k,\varepsilon}$ defined as follows:
$$
g_{k,\varepsilon}(x):=f(x)+\varepsilon\langle h_k^*, x-x_k\rangle\ \ \forall x\in \mathbb{X}.
$$
Then $g_{k,\varepsilon}\in   \mathcal{C}(\mathbb{X}, \mathbb{R})$ and it follows from \eqref{3.24} that $$g_{k,\varepsilon}(z_k)=f(z_k)+\varepsilon\langle h_k^*, z_k-x_k\rangle>0.$$
Thus
$$
0>\inf_{\|h\|=1}d^+g_{k,\varepsilon}(z_k,h)\geq\inf_{\|h\|=1}d^+f(z_k,h)-\varepsilon
=\alpha_k-\varepsilon>-2\varepsilon.
$$
Applying \cref{pro3.1} yields ${\rm Er} g_{k,\varepsilon}\leq 2\varepsilon$, which contradicts (iii), as $\varepsilon>0$ is arbitrary. The proof is complete.\qed

\begin{remark} It should be noted that  the stability of the global error bound for $f$ in \cref{th3.5} may not hold if condition \eqref{3.20} is violated. Consider again the example given in \cref{re3-2}. Then $f(x)=e^x-1$, $S_f=(-\infty, 0]$, ${\rm bdry}(S_f)=\{0\}$ and the stability of the global error bound in \cref{th3.5} is not satisfied. Further, for any $z_k\rightarrow -\infty$, one can verify that
$$
\Big|\inf_{|h|=1}d^+f(z_k,h)\Big|=e^{z_k}\rightarrow 0 \ \ {\rm as} \ k\rightarrow \infty,
$$
which implies that \eqref{3.20} does not hold. \qed
\end{remark}

\subsection{Proofs of \cref{th3.1,th3.2}}

This subsection is devoted to the proofs of \cref{th3.1,th3.2}. We first prove \cref{th3.1}.
\vskip2mm

{\it Proof of \cref{th3.1}} Applying \cite[Theorem 4.2.3]{SIAM13}(or \cite[Proposition 4.5.2]{S2007}), for any $x\in \mathbb{X}$, one has
\begin{equation} \label{4.6}
d^+f(x, h)=\sup_{t\in T_f(x)}d^+f_t(x, h),\ \ \forall h\in \mathbb{X}.
\end{equation}
Denote
$$
\gamma:=\inf\limits_{\|h\|=1}\sup\limits_{t\in T_f(\bar x)}d^+f_t(\bar x, h).
$$

(i)$\Rightarrow$(ii): We first consider the case of $\gamma>0$. By virtue of \cref{pro3.2}, one has  $S_F=\{\bar x\}$. Let $\varepsilon>0$ be such that $\varepsilon< \gamma $ and take $c:= \gamma-\varepsilon>0$. Suppose that $G,g_t$ and $g$ satisfies all conditions  of  (ii). Then for any $t\in T_f(\bar x)\subseteq T_g(\bar x)$, one has
$$
d^+g_t(\bar x,h)\geq d^+f_t(\bar x, h)-\varepsilon.
$$
By using \eqref{4.6} and $T_f(\bar x)\subseteq T_g(\bar x)$, one has
\begin{eqnarray*}
\inf_{\|h\|=1}d^+g(\bar x,h)=\inf_{\|h\|=1}\sup_{t\in T_g(\bar x)}d^+g_t(\bar x,h)&\geq& \inf_{\|h\|=1}\sup_{i\in T_f(\bar x)}d^+g_t(\bar x,h)\\
&\geq& \Big(\inf_{\|h\|=1}\sup_{i\in T_f(\bar x)}d^+f_t(\bar x,h)\Big)-\varepsilon\\
&=&\inf_{\|h\|=1}d^+f(\bar x,h)-\varepsilon=c>0.
\end{eqnarray*}
Using \cref{pro3.2}, one gets
$$
{\rm Er} G(\bar x)={\rm Er} g(\bar x)\geq c.
$$

We next consider the case of $ \gamma<0$.  Let $\varepsilon>0$ be such that $c:=- \gamma -\varepsilon>0$. Then for any $t\in T_g(\bar x)\subseteq T_f(\bar x)$, one has
$$
d^+g_t(\bar x,h)\leq d^+f_t(\bar x, h)+\varepsilon.
$$
This and $T_g(\bar x)\subseteq T_f(\bar x)$ imply
\begin{eqnarray*}
\inf_{\|h\|=1}d^+g(\bar x,h)=\inf_{\|h\|=1}\sup_{t\in T_g(\bar x)}d^+g_t(\bar x,h)&\leq& \inf_{\|h\|=1}\sup_{t\in T_g(\bar x)}(d^+f_t(\bar x,h)+\varepsilon)\\
&\leq& \Big(\inf_{\|h\|=1}\sup_{i\in T_f(\bar x)}d^+f_t(\bar x,h)\Big)+\varepsilon\\
&=&\inf_{\|h\|=1}d^+f(\bar x,h)+\varepsilon=-c<0.
\end{eqnarray*}
Using \cref{pro3.2}, yields
$$
{\rm Er} G(\bar x)={\rm Er} g(\bar x)\geq c.
$$

Note (ii)$\Rightarrow$(iii) follows immediately as $T_f(\bar x)=T_g(\bar x)$. Thus,  it remains to prove (iii)$\Rightarrow$(i).

Let $u^*\in \mathbb{B}_{\mathbb{X}^*}$ and $G_{u^*}\in \mathcal{C}(T\times \mathbb{X}, \mathbb{R})$ be given as in \eqref{4.15}. Note that
$$
g_t(x)=G_{u^*}(t,x)=f_t(x)+\varepsilon\langle u^*, x-\bar x\rangle,
$$
and thus
$$
g(x)=\sup_{t\in T}g_t(x)=\sup_{t\in T}(f_t(x)+\varepsilon\langle u^*, x-\bar x\rangle)=f(x)+\varepsilon\langle u^*, x-\bar x\rangle.
$$
Using the proof in (iii)$\Rightarrow$(i) of \cref{th3.3}, one derives the conclusion. The proof is complete.\qed
\vskip2mm

{\it Proof of \cref{th3.2}} (i)$\Rightarrow$(ii): Based on \cref{pro3.2} and the proof of \cref{th3.1}, we only need to consider the case $\inf_{\|h\|=1}\sup_{t\in T_f(\bar x)}d^+f_t(\bar x,h)\leq 0$ for all $\bar x\in{\rm bdry}(S_F)$.

Using the proof of the claim in \eqref{3.21}, there exists $\varepsilon_0>0$ such that for any $x_0\in{\rm bdry}(S_F)$, one has
\begin{equation}\label{4.24}
\inf\left\{\Big|\inf_{\|h\|=1}d^+f(x,h)\Big|: x\in\mathbb{X}, f(x)\geq-\varepsilon_0\|x-x_0\|\right\}\geq\tau.
\end{equation}
Let $\varepsilon>0$ be such that $\varepsilon<\min\{\varepsilon_0,\tau\}$. Suppose that $G, g_t$ and $g$ satisfy all conditions  as  in (ii). Let $x\in \mathbb{X}$ be such that $g(x)>0$.

We claim that
\begin{equation}\label{4.25}
  \inf_{\|h\|=1}d^+f(x,h)\leq-\tau.
\end{equation}
Granting this, by virtue of $\sup_{t\in T}{\rm Lip}(f_t-g_t)<\varepsilon$ and $T_g(x)\subseteq T_f(x)$, one obtains
\begin{eqnarray*}
\inf_{\|h\|=1}d^+g(x,h)=\inf_{\|h\|=1}\Big(\sup_{t\in T_g(x)}d^+g_t(x,h)\Big)&\leq&\inf_{\|h\|=1}\Big(\sup_{t\in T_g(x)}d^+f_t(x,h)+\varepsilon\Big)\\
&\leq&\inf_{\|h\|=1}\Big(\sup_{t\in T_f(x)}d^+f_t(x,h)\Big)+\varepsilon\\
&=&\inf_{\|h\|=1}d^+f(x,h)+\varepsilon\\
&\leq&-(\tau-\varepsilon).
\end{eqnarray*}
Applying \cref{pro3.2}, yields  ${\rm Er} G={\rm Er}g=\geq \tau-\varepsilon$.

It remains to prove \eqref{4.25}. If $f(x)>0$, then one has \eqref{4.25} by using the proof of \eqref{3.13}. We next consider the case of $f(x)\leq 0$.

By the conditions in (ii), there is $z_0\in{\rm bdry}(S_F)$ such that
$$
f_t(z_0)=g_t(z_0),\ \ \forall t\in T.
$$
Then for any $t\in T_g(x)\subseteq T_f(x)$, one has
\begin{eqnarray*}
f_t(x)\geq g_t(x)-(f_t(z_0)-g_t(z_0))-\varepsilon\|x-z_0\|=g(x)-\varepsilon\|x-z_0\|>-\varepsilon\|x-z_0\|.
\end{eqnarray*}
This implies that $f(x)>-\varepsilon_0\|x-z_0\|$ and thus it follows from \eqref{4.24} that \eqref{4.25} holds.

Note that (ii)$\Rightarrow$(iii) follows immediately since $T_f(x)=T_g(x)$ for all $x\in \mathbb{X}$. It remains to show that (iii) $\Rightarrow$ (i).

Take any $\bar x\in{\rm bdry}(S_F)$ and any $u^*\in \mathbb{B}_{\mathbb{X}^*}$, and let $G\in \mathcal{C}(T\times  \mathbb{X}, \mathbb{R})$ be defined as  in \eqref{4.23}. Note that
$$
g_t(x)=G(t,x)=f_t(x)+\varepsilon\langle u^*, x-\bar x\rangle,
$$
and consequently
$$
g(x)=\sup_{t\in T}g_t(x)=\sup_{t\in T}(f_t(x)+\varepsilon\langle u^*, x-\bar x\rangle)=f(x)+\varepsilon\langle u^*, x-\bar x\rangle.
$$
Thus the conclusion follows from (iii)$\Rightarrow$(i) in the proof of \cref{th3.5}. The proof is complete. \qed


\section{Sensitivity Analysis of Hoffman's Constants for Semi-infinite Linear Systems}

In this section, by using   results from the preceding section, we study sensitivity analysis of Hoffman's constants for semi-infinite linear systems in a Banach space. The aim is to give primal criteria for the Hoffman constants to be uniformly bounded under perturbations on the problem data.

Let $T$ be a compact, possibly infinite metric space (with metric $\rho$) and $a^*:T\rightarrow \mathbb{X}^*, b:T\rightarrow \mathbb{R}$ be continuous functions on $T$. We consider now semi-infinite linear systems in $\mathbb{X}$ defined by
\begin{equation}\label{4-1a}
  \langle a^*(t), x\rangle \leq b(t),\ \ {\rm for\ all} \ t\in T.
\end{equation}
We denote by $\mathcal{S}_{a^*,b}$ the set of solutions to system \eqref{4-1a}. We use the following notations:
$$
\begin{aligned}
&f_{a^*,b}(x):=\max_{t\in T}( \langle a^*(t), x\rangle -b(t)),\\
&J_{a^*,b}(x):=\{t\in T: \langle a^*(t), x\rangle -b(t))=f_{a^*,b}(x)\} \ \  {\rm for \ each } \ x\in \mathbb{X}.
\end{aligned}
$$
It is easy to verify that $J_{a^*,b}(x)$ is a compact subset of $T$ for each $x\in\mathbb{X}$.

Recall that $\mathcal{S}_{a^*,b}$ admits a global error bound, if there exists $\sigma>0$ such that
\begin{equation}\label{4-2a}
 \sigma \mathbf{d}(x, \mathcal{S}_{a^*,b})\leq [f_{a^*,b}(x)]_+, \  {\rm for\ all} \ x\in \mathbb{X}.
\end{equation}
The Hoffman constant of the semi-infinite linear system \eqref{4-1a}, denoted by $\sigma(a^*,b)$, is given by
\begin{equation}\label{4-3a}
\sigma(a^*,b):=\sup\{\sigma>0:  \text{such that } \eqref{4-2a}\ {\rm holds}\}.
\end{equation}
It is easy to verify that
$$
\sigma(a^*,b)=\inf_{x\not\in \mathcal{S}_{a^*,b}}\frac{f_{a^*,b}(x)}{\mathbf{d}(x, \mathcal{S}_{a^*,b})}.
$$

Note that
$$
\mathcal{S}_{a^*,b}={\rm int}(\mathcal{S}_{a^*,b})\cup {\rm bdry}(\mathcal{S}_{a^*,b}).
$$
We begin with characterizing the interior and the boundary of $\mathcal{S}_{a^*,b}$ via the following proposition.

\begin{proposition}\label{pro4.1}
Let be semi-infinite linear system defined by  \eqref{4-1a}. Then
\begin{equation}\label{4-4a}
  {\rm int}(\mathcal{S}_{a^*,b})=\{x\in \mathcal{S}_{a^*,b}: \langle a^*(t), x\rangle < b(t) \ \ {\it for \ all} \ t\in T\}.
\end{equation}
Thus, $x\in {\rm bdry}(\mathcal{S}_{a^*,b})$ if and only if $x\in \mathcal{S}_{a^*,b}$ and there is some nonempty compact  subset $J_x\subseteq T$ such that
\begin{equation}\label{4-5a}
\langle a^*(t), x\rangle = b(t),\forall t\in J_x \ \ {\it and} \ \ \langle a^*(t), x\rangle < b(t),\forall t\in T\backslash J_x.
\end{equation}
\end{proposition}

{\bf Proof.} Let $\bar x\in  {\rm int}(\mathcal{S}_{a^*,b})$. Then there is $\delta>0$ such that $\mathbf{B}(\bar x,\delta)\subseteq  \mathcal{S}_{a^*,b}$.

Suppose on the contrary that there is some $t\in T$ such that
$$
\langle a^*(t), \bar x\rangle = b(t).
$$
Take $u_t\in \mathbb{X}$ such that $\|u_t\|=1$ and $\langle a^*(t), u_t\rangle >0$. Then for any $r\in (0, \delta)$, one has
$$
\langle a^*(t), \bar x+ru_t\rangle >b(t)
$$
which contradicts $\bar x+ru_t\in \mathbf{B}(\bar x,\delta)\subseteq  \mathcal{S}_{a^*,b}$.

On the other hand, let $\bar x\in  \mathcal{S}_{a^*,b}$ be such that $\langle a^*(t), \bar x\rangle < b(t)$ for all $t\in T$. For any $t\in T$, by the continuity of $a^*$ and $b$, there exists $\delta_t>0$ such that
\begin{equation}\label{4.6a}
  \langle a^*(s), x\rangle\leq b(s),\ \ \forall x\in \mathbf{B}(\bar x,\delta_t)\ {\rm and} \ \forall s\in \mathbf{B}_{\rho}(t,\delta_t):=\{t^{\prime}\in T: \rho(t^{\prime},t)<\delta_t\}.
\end{equation}
Using the compactness of $T$, there exist $t_1,\cdots,t_m\in T$ such that
\begin{equation}\label{4.7a}
  T\subseteq\bigcup_{i=1}^m \mathbf{B}_{\rho}(t_i,\frac{\delta_{t_i}}{2}).
\end{equation}
Let $\delta:=\min\{\frac{\delta_1}{2},\cdots,\frac{\delta_n}{2}\}$. We claim that
$$
\mathbf{B}(\bar x,\delta)\subseteq \mathcal{S}_{a^*,b}.
$$
Indeed, for any $x\in\mathbf{B}(\bar x,\delta)$ and any $s\in T$, by \eqref{4.7a}, there exists $j\in\{1,\cdots,m\}$ such that $s\in \mathbb{B}_{\rho}(t_j,\frac{\delta_{t_j}}{2})$. Note that $x\in \mathbf{B}(\bar x,\delta)\subseteq \mathbf{B}(\bar x,\delta_{t_j})$ and it follows from \eqref{4.6a} that
$$
\langle a^*(s), x\rangle\leq b(s).
$$
This implies that the claim holds. The proof is complete. \qed

\vskip 2mm

The following theorem gives primal criteria on the Hoffman's constant for system \eqref{4-1a} under perturbations on the problem data.

\begin{theorem}\label{th4.1}
Consider the following statements:
\begin{itemize}
\item[\rm(a)] Suppose that
\begin{equation}\label{4-6a}
\inf\left\{\Big|\inf_{\|h\|=1}\sup_{t\in J}\langle a^*(t), h\rangle\Big|: J\subseteq T\ {\it is \ compact}\right\}>0.
\end{equation}
Then there exist $c,\varepsilon>0$ such that for all $\widetilde x\in{\rm bdry}(\mathcal{S}_{a^*,b})$ and $\widetilde u^*\in \mathbb{B}_{\mathbb{X}^*}$, one has $\sigma(\widetilde a^*,\widetilde b)\geq c$, where $\widetilde a^*:T\rightarrow \mathbb{X}^*$ and $\widetilde b:T\rightarrow \mathbb{R}$ are defined as follows:
    \begin{equation}\label{4-7a}
      \widetilde a^*(t):=a^*(t)+\varepsilon \widetilde u^*\ \ {\it and} \ \ \widetilde b(t):=b(t)+\varepsilon\langle \widetilde u^*,\widetilde x\rangle\ \ {\it for \ all} \ t\in T.
    \end{equation}
Moreover, the following inequality on Hoffman's constant $\sigma(a^*,b)$ holds:
\begin{equation}\label{4-8a}
\sigma(a^*,b)\geq\inf\left\{\Big|\inf_{\|h\|=1}\sup_{t\in J}\langle a^*(t), h\rangle\Big|: J\subseteq T\ {\it is \ compact}\right\}.
\end{equation}

\item[\rm(b)] Suppose that there exists a sequence $\{x_k\}\subseteq {\rm bdry}(\mathcal{S}_{a^*,b})$ such that
    \begin{equation}\label{4-9a}
\inf_{k\geq 1}\Big|\inf_{\|h\|=1}\sup_{t\in J_{a^*,b}(x_k)}\langle a^*(t), h\rangle\Big|=0.
\end{equation}
Then for any $\varepsilon>0$ there exist $x_{\varepsilon}\in{\rm bdry}(\mathcal{S}_{a^*,b})$ and $u^*_{\varepsilon}\in \mathbb{X}^*$ with $\|u^*_{\varepsilon}\|\leq \varepsilon$ such that $\sigma(a_{\varepsilon}^*,b_{\varepsilon})< \varepsilon$, where $a_{\varepsilon}^*(t):=a^*(t)+u^*_{\varepsilon}$ and $b_{\varepsilon}(t):=b(t)+\langle u^*_{\varepsilon}, x_{\varepsilon}\rangle$ for all $t\in T$.
\end{itemize}
\end{theorem}

{\bf Proof.} We denote
$$
\tau:=\inf\left\{\Big|\inf_{\|h\|=1}\sup_{t\in J}\langle a^*(t), h\rangle\Big|: J\subseteq T\ {\it is \ compact}\right\}>0.
$$
Let $F:T\times\mathbb{X}\rightarrow\mathbb{R}$ be defined by
 $$
 F(t,x):=\langle a^*(t), x\rangle-b(t), \ \forall (t,x)\in T\times\mathbb{X},
$$
and let  $f_t(x):=F(t,x)$ for all $(t,x)\in T\times\mathbb{X}$ and $f$ be given in \eqref{4.3}. Then $f=f_{a^*,b}$, $F\in \mathcal{C}(T\times\mathbb{X},\mathbb{R})$ and each $f_t$ is convex.

$({\rm a})$ We first show that $f_t, f$ and $\tau$ satisfy assumptions as given in \eqref{4.16} and \eqref{3.20a} of \cref{th3.2}.

Indeed, for each  $x\in \mathcal{S}_{a^*,b}$, one has that $J_{a^*,b}(x)$ is a compact subset of $T$ and thus
$$
\inf\left\{\Big|\inf_{\|h\|=1}\sup_{t\in J_{a^*,b}(x)}\langle a^*(t), h\rangle\Big|: x\in \mathcal{S}_{a^*,b} \right\}\geq\tau.
$$
Note that
$$
\inf_{\|h\|=1}d^+f(x,h)=\inf_{\|h\|=1}\sup_{t\in T_f(x)}d^+f_t(x,h)=\inf_{\|h\|=1}\sup_{t\in J_{a^*,b}(x)}\langle a^*(t), h\rangle
$$
and it follows from \eqref{4-6a} that $f_t, f$ and $\tau$ satisfy assumptions in \eqref{4.16} and \eqref{3.20a}.

Based on the equivalence of (i) and (iii) in \cref{th3.2}, one derives the conclusion.
\vskip 2mm
We next prove \eqref{4-8a}. Using the proof of (i)$\Rightarrow$(ii) in \cref{th3.2}, for any $\varepsilon\in (0,\tau)$ sufficiently small, one has (ii)  as given in \cref{th3.2} holds for $c:=\tau-\varepsilon$, and this implies that
$$
\sigma(a^*,b)\geq c=\tau-\varepsilon.
$$
Letting $\varepsilon\rightarrow 0^+$,  one obtains \eqref{4-8a}.

$({\rm b})$  Suppose that there exists a sequence $\{x_k\}\subseteq {\rm bdry}(\mathcal{S}_{a^*,b})$ such that \eqref{4-9a} holds. Then using the proof of (iii)$\Rightarrow$(i) in \cref{th3.2}, yields the conclusion. The proof is complete.\qed

\vskip 2mm

As an application of \cref{th4.1}, let us give an example where one can verify the uniform bound of the Hoffman constant for the given linear system.

\begin{example}
Let $\mathbb{X}:=\mathbb{R}^2$, $T:=\{1,2,3\}$, let $a^*:T\rightarrow \mathbb{R}^2, b:T\rightarrow \mathbb{R}$ be defined as follows:
$$
a^*(1):=(1,1), a^*(2):=(-2,1), a^*(3):=(1,-2), b(1):=1, b(2):=2, b(3):=2.
$$
Consider the following linear system:
\begin{equation}\label{4-10a}
  \langle a^*(t), x\rangle \leq b(t),\ \ t\in T.
\end{equation}
Take $J_1:=\{1,2\}, J_2:=\{2,3\}, J_3:=\{1,3\}, J_4:=\{1,2,3\}$. Consider the following four min-max optimization problems:
\begin{equation}
OP(J_k) \ \ \ \
 \min_{\|h\|=1} \max_{t\in J_k}\langle a^*(t),h\rangle, k=1,2,3,4.
\end{equation}

We denote by $\theta_k$ and $h_k$ the optimal value and the optimal solution of $OP(J_k)$ for each $k$, respectively. Solving  these four problems we obtain that
$$
\begin{aligned}
&\theta_1=-1, h_1=(0, -1); \theta_2=-\frac{\sqrt{2}}{2}, h_2=(-\frac{\sqrt{2}}{2}, -\frac{\sqrt{2}}{2});\\ &\theta_3=-1, h_3=(-1, 0); \theta_4=\sqrt{2}, h_4=(-\frac{\sqrt{2}}{2}, -\frac{\sqrt{2}}{2}).
\end{aligned}
$$
This implies that
$$
\min\left\{\Big|\min_{\|h\|=1}\max_{t\in J}\langle a^*(t), h\rangle\Big|: J\subseteq T\right\}=\frac{\sqrt{2}}{2}>0,
$$
and then \cref{th4.1} implies that the Hoffman constant ($\sigma(\cdot)$) for system \eqref{4-10a} is uniformly bounded under perturbations on the problem data $(a^*,b)$.

\end{example}

Finally, let us  show  that the uniform bound of the Hoffman constant may not be satisfied if the condition \eqref{4-6a} is violated in  \cref{th4.1}.

\begin{example}
Let $\mathbb{X}:=\mathbb{R}^2$, $T:=\{1,2\}$, and let $a^*:T\rightarrow \mathbb{R}^2, b:T\rightarrow \mathbb{R}$ be defined as follows:
$$
a^*(1):=(1,1), a^*(2):=(-1,-1), b(1)=b(2):=0.
$$
Consider the following linear system:
\begin{equation}\label{4-12a}
  \langle a^*(t), x\rangle \leq b(t),\ \ t\in T.
\end{equation}
Then $\mathcal{S}(a^*,b)=\{x=(x_1, x_2)\in \mathbb{R}^2: x_1+x_2=0\}$. Note that
\begin{equation*}
 \min_{\|h\|=1} \max_{t\in T}\langle a^*(t),h\rangle=0
\end{equation*}
and thus from \cref{th4.1}, the Hoffman constant for system \eqref{4-12a} is not uniformly bounded under perturbations on the problem data $(a^*,b)$.

Indeed, for any $\varepsilon>0$,   select $u_{\varepsilon}^*:=(0,\varepsilon)$, $x_{\varepsilon}:=(0,0)$ and define $a_{\varepsilon}^*(t):=a^*(t)+u_{\varepsilon}^*,b_{\varepsilon}(t):=b(t)$ for $t=1,2$. Then $f_{a_{\varepsilon}^*,b_{\varepsilon}}(x):=\max\{x_1+x_2+\varepsilon x_2, -x_1-x_2+\varepsilon x_2\}$ for all $x=(x_1, x_2)\in \mathbb{R}^2$ and
$$
\mathcal{S}(a_{\varepsilon}^*,b_{\varepsilon})=\{x=(x_1, x_2)\in \mathbb{R}^2: x_1+(1+\varepsilon)x_2\leq 0, x_1+(1-\varepsilon)x_2\geq 0\}.
$$
Take $u_{\varepsilon}:=(-\varepsilon, \varepsilon)$ and thus
$$
\frac{f_{a_{\varepsilon}^*,b_{\varepsilon}}(u_{\varepsilon})}{\mathbf{d}(u_{\varepsilon}, \mathcal{S}(a_{\varepsilon}^*,b_{\varepsilon}))}=\frac{\varepsilon^2}{\sqrt{2}\varepsilon}
=\frac{\varepsilon}{\sqrt{2}},
$$
which implies that  $\sigma(a_{\varepsilon}^*,b_{\varepsilon})<\varepsilon$.
\end{example}

\section{Concluding Remarks and Perspectives}

This paper is devoted to studying stability of local and global error bounds for semi-infinite convex constraint systems in a Banach space. We establish primal characterizations of the stability of error bounds in terms of directional derivatives. It is proved that verifying the stability of the error bounds for semi-infinite convex constraint systems is equivalent to solving several minimax problems defined by directional derivatives of component functions. Stability under perturbations of the problem data only requires that all component functions have the same linear perturbations. When applied to sensitivity analysis of Hoffman's constants for semi-infinite linear systems, these stability results lead to conditions ensuring that Hoffman's constants can be uniformly bounded.

In a future work, we plan to study applications of these results to sensitivity analysis of variational problems and convergence analysis of computational algorithms. Primal criteria of stability of error bounds for inequality systems defined by non-convex functions are also to be considered.

%

\bibliography{WTY}

\bibliographystyle{unsrt}

\end{document}